\documentclass[letterpaper,12pt]{article}
\usepackage{amssymb}
\usepackage{amsmath}
\usepackage{amscd}
\usepackage[all]{xy}

\newtheorem{theo}{Theorem}[section]
\newtheorem{prop}[theo]{Proposition}
\newtheorem{lem}[theo]{Lemma}
\newtheorem{cor}[theo]{Corollary}
\newtheorem{defi}[theo]{Definition}

\def \Br {{\rm{Br}}}

\def \Ga {{\Gamma}}
\def \R {{\mathbb{R}}}

\def \val {{\rm {val}}}

\def \Gal {{\rm{Gal}}}
\def \Ker {{\rm{Ker\,}}}

\def \A{{\mathbb A}}
\def \P{{\mathbb P}}

\def \dim {{\rm{dim\,}}}
\def \ev {{\rm{ev}}}
\def \Hom {{\rm {Hom}}}
\def \End {{\rm {End}}}

\def \GL {{\rm {GL}}}

\def\ov{\overline}

\def \Z {{\mathbb Z}}
\def \Q {{\mathbb Q}}
\def \F {{\mathbb F}}

\def\T{{\cal T}}

\def\lra{\longrightarrow}

\def\H{{\rm H}}

\def\O{{\cal O}}

\def\inv{{\rm inv}}

\def\Kum{{\rm Kum}}

\def\O{{\cal O}}

\def\sA{{\mathcal A}}

\newcommand{\bthe}{\begin{theo}}
\newcommand{\ble}{\begin{lem}}
\newcommand{\bpr}{\begin{prop}}
\newcommand{\bco}{\begin{cor}}
\newcommand{\bde}{\begin{defi}}
\newcommand{\ethe}{\end{theo}}
\newcommand{\ele}{\end{lem}}
\newcommand{\epr}{\end{prop}}
\newcommand{\eco}{\end{cor}}
\newcommand{\ede}{\end{defi}}

\def\sB{{\mathcal B}}
\def\sC{{\mathcal C}}

\def\beq{\begin{equation} \label}

\title{Odd order Brauer--Manin obstruction on diagonal quartic surfaces}
\author{Evis Ieronymou and Alexei N. Skorobogatov}
\date{}
\begin{document}
\baselineskip=15pt
\maketitle

\begin{abstract}
We determine the odd order torsion subgroup of the Brauer
group of diagonal quartic surfaces over the field of rational
numbers. We show that a non-constant Brauer element of 
odd order always obstructs weak
approximation but never the Hasse principle.
\end{abstract}

\section{Introduction}

This paper is devoted to the arithmetic of diagonal
quartic surfaces over the field of rational numbers. Let $D=[a,b,c,d]\subset\P^3_\Q$
be the surface given by
\begin{equation}
ax^4+by^4=cz^4+dw^4, \label{*}
\end{equation}
where $a,\,b,\,c,\,d\in\Q^*$. The set of rational points
$D(\Q)$ is a subset of the set of adelic points 
$D(\A_\Q)=\prod D(\Q_p)$, where the product is over all
completions of $\Q$ including $\R$.
The closure of $D(\Q)$ in the product topology
is contained in the Brauer--Manin set $D(\A_\Q)^\Br$ 
defined as the set of adelic
points orthogonal to the Brauer group $\Br(D)$ with respect to the pairing
provided by class field theory, see \cite[Ch. 6]{Sk}.
The image $\Br_0(D)$ of the natural map $\Br(\Q)\to\Br(D)$ pairs trivially
with $D(\A_\Q)$. The group $\Br(D)/\Br_0(D)$ is finite \cite{SZ1}.
However, the computation of
this group for an arbitrary diagonal quartic $D$ over $\Q$ is still an
open problem. Results obtained in \cite{Evis, ISZ} allow one to prove that
$\Br(D)/\Br_0(D)$ is zero in certain particular cases. 
Recall that $\Br_1(D)$ denotes the kernel of the natural map
$\Br(D)\to \Br(\ov D)$. M. Bright \cite{Bright, B06}
determined the group structure of $\Br_1(D)/\Br_0(D)$
for any diagonal quartic surface $D$ over $\Q$,
in particular, he showed that the order of this group is a power of 2.

In this paper we compute
the odd order torsion subgroup $\big(\Br(D)/\Br_0(D)\big)_{\rm odd}$ 
of $\Br(D)/\Br_0(D)$,
improving on the estimate obtained previously by the authors 
jointly with Yuri Zarhin \cite[Cor. 3.3, Cor. 4.6]{ISZ}. 
Let $\Ga_\Q={\rm Gal}(\ov \Q/\Q)$ be the Galois group of $\Q$,
and let $\ov D=D\times_\Q\ov\Q$. Let
$\langle a\rangle\subset\Q^*$ be the cyclic subgroup generated by $a\in\Q^*$.

\bthe \label{main}
Let $D=[a,b,c,d]$ be a diagonal quartic surface over $\Q$. Then
$$\big(\Br(D)/\Br_0(D)\big)_{\rm odd}=\Br(\ov D)_{\rm odd}^{\Ga_{\Q}}\simeq
\left \{ \begin{array}{ll}
\Z/3 & \textrm{if}  \ -3abcd \in \langle-4\rangle\Q^{*4},\\
\Z/5 & \textrm{if}  \ 125abcd \in \langle-4\rangle\Q^{*4},  \\
0 & \textrm{otherwise.}
\end{array}
\right .
$$
\ethe

For an odd prime $\ell$ the results of M. Bright and the
divisibility of every element of $\Br_0(D)$ by $\ell$ imply
the surjectivity of the natural map
$\Br(D)_\ell\to(\Br(D)/\Br_0(D))_\ell$. Thus for $\ell=3$
or $\ell=5$ under the conditions of Theorem \ref{main}
we have an element of $\Br(D)$ of order $\ell$ which is not
in $\Br_0(D)$. The first example of a diagonal quartic surface
$D$ with a non-constant 3-torsion element in $\Br(D)$ 
was found by Thomas Preu \cite{Preu}. 
Preu's surface $D=[1,3,4,9]$ has an obvious rational point, 
and he shows that for any rational point on his surface one has 
$3|yw$ while there are $\Q_3$-points
not satisfying this condition. This failure of weak approximation is
explained by a 3-torsion element in $\Br(D)$. We show that
this situation is quite general.

Let $\inv_p:\Br(\Q_p)\tilde\lra\Q/\Z$ be
the local invariant isomorphism of class field theory.
For an element $\sA\in\Br(D)_n$ 
we denote by $\ev_{\sA,p}:D(\Q_p)\to \frac{1}{n}\Z/\Z$ 
the evaluation map at the prime $p$, defined as
$\ev_{\sA,p}(P)=\inv_p(\sA(P))$.

\bthe \label{main2}
Let $D$ be a diagonal quartic surface over $\Q$ such that 
$D(\Q_p)\not=\emptyset$. If $\sA\in\Br(D)_\ell\setminus\Br_0(D)$, 
where $\ell$ is an odd prime, then

{\rm (i)} the map $\ev_{\sA,p}$ is constant for $p\neq \ell$,

{\rm (ii)} the map $\ev_{\sA,\ell}$ is surjective.
\ethe
{\em Proof.} By Theorem \ref{main} we only need to consider
$\ell=3$ and $\ell=5$. The statement (i) is proved
in Proposition \ref{Pr5fo}, and the statement (ii) 
is proved in Proposition \ref{Pr5f5} for $\ell=5$ and in
Proposition \ref{Pr3f3} for $\ell=3$. 
$\Box$

\medskip

Since $\Br(\R)\cong\Z/2$, evaluating an element of 
$\Br(D)_{\rm odd}$ at any point of $X(\R)$ gives zero,
hence we obtain the following result.

\bco \label{main3}
If $D$ is a diagonal quartic surface over $\Q$
such that $D(\A_\Q)$ is non-empty, then 
$D(\A_\Q)^{\Br(D)_{\rm odd}}$ is non-empty too.
If $\Br(D)_{\rm odd}$ is not contained in $\Br_0(D)$, then
$D(\A_\Q)^{\Br(D)_{\rm odd}}\not=D(\A_\Q)$.
\eco

The fact that $\Br(D)_{\rm odd}$ never obstructs the Hasse principle
and always obstructs weak approximation on diagonal quartic surfaces
over $\Q$ came to us as a surprise.  
It seems that all known counterexamples to the Hasse principle
on K3 surfaces are given by elements of even order, see 
\cite{Bright, B06, HV} and references in these papers. 
This prompts the following general question. 
Does there exist a K3 surface $X$ over a number field $k$ with $X(\A_k)\not=\emptyset$
such that $X(\A_k)^{\Br(X)_{\rm odd}}$ is empty? 
Note that when the degree of the polarisation of $X$ 
(which is always even) is a power of 2, e.g. when $X$ is a quartic
K3 surface, there is a 0-cycle of degree 1 on $X$ over each
completion of $k$ such that this collection is orthogonal to 
$\Br(D)_{\rm odd}$ with respect to the Brauer--Manin pairing.
For everywhere locally soluble diagonal quartic surfaces $D$ over $\Q$ 
whose coefficients $a,b,c,d$ are general enough, M. Bright \cite{B11}
has shown using \cite{Evis, ISZ}
that $\Br(D)/\Br_0(D)=\Br_1(D)/\Br_0(D)\simeq\Z/2$ but the Brauer
group does not obstruct the Hasse principle. Conditional results
on the existence of rational points on diagonal quartic surfaces
can be found in \cite{SD} and \cite[Thm. 1.51]{W}.

The failure of weak approximation in Corollary \ref{main3}
can be illustrated in explicit examples.
The following result is an immediate corollary of the proof 
of Proposition~\ref{Pr5f5}.

\bco \label{co1}
Let $D=[a,b,c,d]$, where $a,b,c,d$ are fourth power free integers
such that $5^3abcd\in\langle -4\rangle\Q^{*4}$.
Let $x,y,z,w$ be integers such that $(x,y,z,w)\in D(\Q)$.

{\rm (1)} If $a,b,c\in\Z_5^*$ and $a+b\equiv 0\bmod 5$, then $25|xyzw$.

{\rm (2)} If $a,b\in\Z_5^*$ and $c\equiv d\equiv 0\bmod 5$, then $5|zw$.
\eco

We would like to emphasise that these divisibility conditions cannot be obtained
by solving corresponding congruences modulo powers of $5$. 
Indeed, any surface satisfying the assumption of (1) 
has a solution $(\alpha,1,5,1)$ for some $\alpha\in\Z_5^*$, 
and any surface satisfying the assumption of 
(2) has a solution $(\beta,1,1,1)$ for some $\beta\in\Z_5^*$.
For example, for any $n\in\Z$ coprime to $5$ and any $\epsilon\in\{0,1\}$
the surfaces in the family
$$x^4-y^4=n z^4-5n^3(-4)^{\epsilon}w^4$$
satisfy the condition in (1), hence all integral solutions
of these equations satisfy $25|xyzw$.
The surfaces 
$$x^4-y^4= 5^2n z^4-5^3n^3(-4)^{\epsilon}w^4$$
satisfy the condition in (2), so in this case we have $5|zw$. 

Similarly, the next statement is a corollary of the proof 
of Proposition \ref{Pr3f3}.

\bco \label{co3}
Let $D=[a,b,c,d]$, where $a,b,c,d$ are fourth power free integers
such that $-3abcd\in\langle -4\rangle\Q^{*4}$.
Let $x,y,z,w$ be integers such that $(x,y,z,w)\in D(\Q)$.

{\rm (1)} If $a,b,c\in\Z_3^*$ and $a\equiv b\equiv c\bmod 3$, then $9|xyw$.

{\rm (2)} If $a,b\in\Z_3^*$ and $c\equiv d\equiv 0\bmod 3$, then $3|zw$.
\eco

In the same way as before
these divisibility conditions cannot be obtained
by solving corresponding congruences
modulo powers of $3$. For any $n\equiv 1\bmod 3$ 
and any $\epsilon\in\{0,1\}$ the surfaces in the family
$$x^4+ny^4=z^4-27(-4)^{\epsilon}n^3w^4$$
satisfy the condition in (1), hence all integral solutions
of these equations satisfy $9|xyw$. (For any $n\equiv 2\bmod 3$
we have $9|yzw$.) The surfaces 
$$x^4-y^4=3n z^4+9(-4)^{\epsilon}n^3w^4$$
satisfy the condition in (2), so in this case we have $3|zw$. 
The same holds for the surface $x^4-4y^4=-3 z^4+9w^4$ considered 
in \cite{Preu}.

The proof of Theorem \ref{main} uses a geometric step which is
a reduction to a certain Kummer surface and then to
a product of two elliptic curves $E^{m_1}\times E^{m_2}$, see
Section~\ref{section2}. 
These curves have complex multiplication
by $\sqrt{-1}$ and the Weierstrass equation $y^2=x^3-mx$, so they are
quartic twists of each other. The Galois module structure of 
$\Br(E^{m_1}\times E^{m_2})$ is then explicitly described in 
Section \ref{section3} leading to a proof of Theorem \ref{main}
in Section \ref{section4}. The proof of Theorem \ref{main2}
is carried out in Section~\ref{section5}. Our main tool is the standard
theory of the formal group attached to an elliptic curve; it
provides a simple way to decide if a point on an elliptic curve
over $\Q_p$ is divisible by $\ell$. The case when $p\not=\ell$
is easy, but the case $p=\ell$, although it is essentially a finite
calculation, requires a subtle analysis of some local pairings.

We expect the methods of this paper to be applicable to more
general surfaces dominated by
products of curves, such as the surfaces in $\P^3$
given by an equation $P(x,y)=Q(z,w)$, where $P$ and $Q$
are homogeneous polynomials of the same degree. 

This paper started as an attempt to
understand the calculations in \cite{Preu}, and we are grateful to
Thomas Preu for a very useful explanation of his work at
the conference ``Arithmetic of surfaces"
at the Lorentz Center (Leiden). The work on this paper continued during the
second author's visits to
the University of Cyprus (Nicosia), the Centre Interfacultaire
Bernoulli (Lausanne) and the Hausdorff Research Institute for Mathematics
(Bonn). He is deeply grateful to these institutions for their hospitality.
We thank Rachel Newton for her interest in this paper and
Anthony V\'arilly-Alvarado for useful discussions.

\section{Isogenies between K3 and abelian surfaces} \label{section2}

Let $k$ be a field of characteristic $0$ with an algebraic closure
$\ov k$ and absolute Galois group $\Ga_k=\Gal(\ov k/k)$.
For a variety $X$ over $k$ we write $\ov X=X\times_k\ov k$.

A rational dominant map of K3 or abelian surfaces
is called an {\it isogeny}. It is well known that the ranks of
the lattices of transcendental cycles of isogenous surfaces are equal.
This implies that the geometric Brauer group of a K3 or abelian surface,
as an abelian group, depends only on the isogeny class of the surface. 
Indeed, it is isomorphic to $(\Q/\Z)^t$,
where $t$ is the rank of the lattice of transcendental cycles.

We now consider isogenies between certain
quartic K3 surfaces in $\P^3$ and products of two
curves of genus 1.
Let $p(t)$ and $q(t)$ be separable polynomials of degree
$4$. We denote the genus 1 curves
$$y^2=p(x), \quad \quad u^2=q(v)$$ by $C$ and $C'$, respectively.
Let $J={\rm Jac}(C)$ and $J'={\rm Jac}(C')$.
It is a classical fact that the Jacobian $J$ of $C$ is given by $u^2=f(t)$,
where $f(t)$ is the (monic) resolvent cubic polynomial of $p(x)$,
see \cite[Lemma 3]{BSD} or \cite[Prop. 3.3.6 (a)]{Sk}.

Let $P(x,y)$ and $Q(z,u)$ be the quartic forms such that
$p(t)=P(t,1)$ and $q(t)=Q(t,1)$, and 
let $X$ be the quartic surface in $\P^3$ given by 
$$P(x,y)=Q(z,w).$$
It is smooth, and is therefore a K3 surface.
Let $Y$ be the Kummer surface $\Kum(C\times C')$ defined as the minimal desingularisation
of the quotient of $C\times C'$ by the involution simultaneously
changing the signs of $y_1$ and $y_2$.
Finally, let $Z=\Kum(J\times J')$.

\bpr \label{p1}
There is a commutative diagram of isogenies whose degrees
are powers of $2$:
\begin{equation}
\xymatrix{&C\times C' \ar[r] \ar@{-->}[d]& J\times J'\ar@{-->}[d]\\
X \ar@{-->}[r]  &Y\ar@{-->}[r]   &Z
}\label{d1}
\end{equation}
\epr
{\it Proof} The vertical maps are isogenies of
degree 2 that come from the definitions of $Y$ and $Z$.
The isogeny $\xymatrix{X \ar@{-->}[r]&Y}$ of degree 2 is constructed
as follows. The surface
$Y$ is birationally equivalent to the affine surface $u^2p(t)=q(z)$.
Setting $u=y^2$, $yt=x$, we obtain $P(x,y)=Q(z,1)$, which is an
affine equation of $X$.

Choosing a root of $p(x)=0$ as the origin of the group law on $\ov C$
we define an
isomorphism $\ov C\cong\ov J$ that identifies the action of $J$ on $C$
with the action of $J$ on $J$ by translations. This isomorphism
identifies the involution $x\mapsto -x$ on $J$ with the involution
on $C$ that changes the sign of $y$. The action of $J_2$ on $C$
gives rise to the \'etale morphism $C\to J=C/J_2$. It is compatible
with the involutions on $C$ and $J$. Thus we have a finite \'etale morphism
$C\times C'\to J\times J'$ of degree 16, which is a torsor with 
structure group $J_2\times J'_2$. Since this morphism is compatible 
with the involutions, it induces an isogeny $\xymatrix{Y\ar@{-->}[r]& Z}$
of degree 16. $\Box$

\bco \label{c1}
Diagram {\rm (\ref{d1})} induces the following
commutative diagram of homomorphisms of $\Ga_k$-modules
\begin{equation}
\xymatrix{\Br(\ov J\times \ov J')\ar[r]&
\Br(\ov C\times \ov C')&\\
\Br(\ov Z)\ar[r]\ar[u]_{\cong}&\Br(\ov Y)\ar[r]\ar[u]_{\cong}&\Br(\ov X)}\label{d2}
\end{equation}
All the homomorphisms in this diagram induce bijections
on the subgroups of elements of odd order.
\eco
{\em Proof.} In \cite[Prop. 1.1]{ISZ} it is proved that
if $V$ and $W$ are geometrically integral smooth varieties over
$k$, and $f:V\to W$ is a dominant, generically
finite morphism of degree $d$, then the kernel of
$f^*:\Br(W)\to\Br(V)$ is killed by $d$.
If $f:\xymatrix{X\ar@{-->}[r]& Y}$ is an isogeny of K3 or abelian 
surfaces of degree $d$,
then there is a birational morphism $\varphi:X'\to X$ and a finite
surjective morphism $f':X'\to Y$ of degree $d$ such that $f=f'\varphi^{-1}$.
Since $\varphi^*:\Br(\ov X)\to\Br(\ov X')$ is an isomorphism,
there is an induced map of $\Ga_k$-modules
$f^*:\Br(\ov Y)\to\Br(\ov X)$. This map is surjective with the kernel
annihilated by $d$.
In the particular case when $A$ is an abelian surface,
the natural isogeny
$f:\xymatrix{A\ar@{-->}[r]& Y}$ induces an isomorphism $f^*:\Br(\ov Y)\to
\Br(\ov A)$, by \cite[Prop. 1.3]{SZ2}.
Thus the vertical arrows are isomorphisms. The last statement follows
from the fact that the degrees of all the isogenies in (\ref{d1}) are powers of 2.
$\Box$


\bco \label{c2}
Let $n$ be an odd positive integer. Then the $\Ga_k$-module $\Br(\ov X)_n$
is canonically isomorphic to the $\Ga_k$-module
$$\Br(\ov J\times \ov J')_n=\Hom(J_n,J'_n)/\big(\Hom(\ov J,\ov J')/n\big).$$
\eco
{\it Proof} This follows from the last sentence of Corollary \ref{c1}
and \cite[Prop. 3.3]{SZ2}
which gives the equality in the displayed formula. $\Box$

\medskip

\noindent{\bf Remark} Examples of elliptic curves $J$, $J'$
over $\Q$ such that
$\Br(\ov J\times \ov J')^{\Ga_\Q}_{\rm odd}=0$ were
constructed in \cite[Section 4]{SZ2}. By
Corollary \ref{c2} these examples give rise to
quartic surfaces $X\subset\P^3_\Q$
with $\Br(X)_{\rm odd}\subset\Br_0(X)$.

\section{Twists of lemniscata} \label{section3}

Let us denote by $E^c$ the elliptic curve $y^2=x^3-cx$, where $c\in k^*$.

\bco\label{mco}
Let $D=[a,b,c,d]$ be a diagonal quartic surface over $k$.
If $n$ is odd, then the $\Ga_k$-module $\Br(\ov D)_n$
is canonically isomorphic to
$$\Hom(E^{4ab}_n,E^{4cd}_n)/\big(\Hom(\ov{E^{4ab}} ,\ov{E^{4cd}})/n\big).$$
\eco
{\em Proof.} In our previous notation we have $p(t)=at^4+b$ and $q(t)=ct^4+d$.
From \cite[Lemma 3]{BSD} or \cite[Prop. 3.3.6 (a)]{Sk} we see
that $J$ is given by $y^2=x^3-4abx$ and $J'$ is given by 
$y^2=x^3-4cdx$, so that $J=E^{4ab}$, $J'=E^{4cd}$.
It remains to apply Corollary \ref{c2}. $\Box$

\medskip

The curve $E^c$ has complex multiplication by $\O=\Z[\sqrt{-1}]$, 
so we need to
consider in detail elliptic curves with such complex multiplication.
The {\em lemniscata} is the elliptic curve $E$ with the equation $y^2=x^3-x$.
Note that $\O$ is the full endomorphism ring $\End(E)$ of $E$.
The normalised action of $\O$ on $E$ is the action
such that the induced action of $z\in\O$ on the differential form
$y^{-1}dx$ is the multiplication by $z$. Then $i=\sqrt{-1}$
sends $(x,y)$ to $(-x,iy)$. For $z\in \O$ we
denote the normalised action of $z$ by $[z]\in\End(E^c)$.

Since $\O^*=\{\pm 1,\pm i\}=\mu_4$, we see that
the group $k$-scheme $\mu_4$ acts on $E$. The following well known lemma
is checked directly from the definition of twisting.

\ble \label{le1}
The curve $E^c$ is the quartic twist of $E$ corresponding to the
class of $c^{-1}$ in $\H^1(k,\mu_4)=k^*/k^{*4}$.
\ele

Let $n$ be an odd positive integer.
If $\phi:E_n\to E_n$ is an endomorphism, and $x,\,y\in\O^*$, then
$y\phi x^{-1}$ is also an endomorphism of $E_n$. Thus
$\End(E_n)$ has a natural structure of a $\mu_4\times \mu_4$-module.
Consider the action of $[i]$ on $\End(E_n)$ by conjugation, that is,
via the diagonal map $\mu_4\to \mu_4 \times \mu_4$ .
Since all endomorphisms commute with $[-1]$,
and $n$ is odd, $\End(E_n)$ is the
direct sum of the eigenspaces of $[i]$ with eigenvalues
$1$ and $-1$:
\begin{equation}
\End(E_n)=\End(E_n)^+\oplus \End(E_n)^-. \label{e1}
\end{equation}
This is clearly a direct sum of $\mu_4\times\mu_4$-modules.
It is easy to see that ${\rm rk}(\End(E_n)^+)={\rm rk}(\End(E_n)^-)=2$.
It follows that the inclusion
$\O/n\subset\End(E_n)^+$ is an equality $\End(E_n)^+=\O/n$.
The action of the Galois group $\Ga_k$ preserves this
decomposition, hence we have canonical isomorphisms
of $\Ga_k$-modules, and also of $\mu_4\times\mu_4$-modules:
\begin{equation}
\End(E_n)^-=\End(E_n)/(\O/n).\label{i1}
\end{equation}
By Corollary \ref{mco} this is canonically isomorphic to the $\Ga_k$-module
$\Br(\ov X)_n$, where $X=[1,1,1,1]$ is the `untwisted' diagonal quartic surface.
We note that the actions of $\mu_4$ on $\End(E_n)^-$ coming from the action
of $\mu_4$ on the source $E_n$ and on the target $E_n$, coincide.
(Indeed, $[i]\varphi=\varphi[i]^{-1}$ for any $\varphi\in \End(E_n)^-$.)
We shall always consider $\End(E_n)^-$ with this structure of a $\mu_4$-module.

\bpr\label{p2}
Let $D=[a,b,c,d]$ be a diagonal quartic surface over $k$.
If $n$ is odd, then the $\Ga_k$-module $\Br(\ov D)_n$
is canonically isomorphic to the twist of $\End(E_n)^-$
by a cocycle whose class in $\H^1(k,\mu_4)=k^*/k^{*4}$
is represented by $(abcd)^{-1}$.
\epr
{\em Proof.} We use Corollary \ref{mco}. By Lemma \ref{le1} the
$\Ga_k$-module $\Hom(E^{4ab}_n,E^{4cd}_n)$ is the twist
of $\End(E_n)$ by a cocycle with the class
$(4ab,4cd)^{-1}\in \H^1(k,\mu_4)\times \H^1(k,\mu_4)$,
with respect to its natural $\mu_4\times\mu_4$-module structure.
Now the statement follows from the isomorphism (\ref{i1}),
since the right and left actions of $\mu_4$ on $\End(E_n)^-$ coincide. $\Box$

\medskip

\noindent{\bf Remark} Any quartic twist of $E$ has complex multiplication
by $\O$, so similarly to (\ref{e1}) we can consider a decomposition
of $\Ga_k$-modules
$$
\Hom(E_n^{ab},E_n^{cd})=\Hom(E_n^{ab},E_n^{cd})^+\oplus
\Hom(E_n^{ab},E_n^{cd})^-.
$$
Then Proposition \ref{p2} and Corollary \ref{c1} give canonical isomorphisms
of $\Ga_k$-modules
\begin{equation}
\Br(\ov D)_n=\Br(\ov Z)_n=\Hom(E_n^{ab},E_n^{cd})^-=\Hom(E_n^{cd},E_n^{ab})^-.
\label{e2}
\end{equation}

\bpr \label{c3}
Let $D=[a,b,c,d]$ be a diagonal quartic surface over $k$.
If $n$ is odd, then $\big(\Br(D)/\Br_0(D)\big)_n=\Br(\ov D)^{\Ga_k}_n.$
\epr
{\em Proof.} Consider the natural map
$$\Br(D)/\Br_0(D) \to \Br(\ov D)^{\Ga_k}.$$
According to Example 1 after Prop. 5.1 in
\cite[Thm. 4.3]{CS} the cokernel of this map is a finite abelian 2-group.
On the other hand, the calculations in M. Bright's thesis \cite{Bright}
show that the kernel is also a finite abelian 2-group. $\Box$

\section{Proof of Theorem \ref{main}} \label{section4}

It is known that in the proof of Theorem \ref{main} we only need to consider
torsion subgroups of order 3 and 5. Indeed, we have the following

\ble\label{pp}
For any odd prime $\ell $ we have
$\Br(\ov D)^{\Ga_\Q}_{\ell^\infty}=\Br(\ov{D})^{\Ga_\Q}_\ell$.
If $\ell\geq 7$, then $\Br(\ov{D})^{\Ga_\Q}_\ell=0$.
\ele
{\em Proof.} \cite[Prop. 4.1, Thm. 3.2]{ISZ}. $\Box$

\ble \label{lem3} For $\ell=3$ or $\ell=5$
consider $\O/\ell$ with its natural structure of a $\Ga_\Q$-module and the
$\Ga_\Q$-compatible action of $\O^*=\mu_4$. Then $\End(E_3)^-$
is isomorphic to the quartic twist of $\O/3$ by a cocycle with the class
$-1/3\in \H^1(\Q,\mu_4)=\Q^*/\Q^{*4}$, and $\End(E_5)^-$
is isomorphic to the quartic twist of $\O/5$ by a cocycle with the class
$5\in \Q^*/\Q^{*4}$. Moreover these are also $\mu_4$-module 
isomorphisms, with respect to the action of the twisted group 
(which is canonically $\mu_4$) on the twisted module.
\ele
{\em Proof.}
Recall that an irreducible element $\pi\in\O$ is called a {\em primary prime}
if $\pi\equiv 1 \bmod (1+i)^3$. For a primary prime
$\pi$ let ${\rm Frob}_\pi\in \Gal(\Q(i)^{\rm{ab}}/\Q(i))$ be the image 
under the Artin map
of the id\`{e}le with component $\pi$ at the prime ideal generated by 
$\pi$ and component $1$ at all other places of $\O$.
The Hecke character associated to $E/\Q(i)$ takes the value $\pi$ at 
the prime ideal generated by $\pi$ (\cite[Exercise 2.34]{S2})
Hence we get the classical fact (Gauss, Deuring) that the action of 
$\Ga_{\Q(i)}$ on $E_{\ell}$ factors through the homomorphism
$\Ga_{\Q(i)}\to (\O/\ell)^*$ that sends ${\rm Frob}_\pi$
to the class of $\pi$ in $(\O/\ell)^*$, where $\pi$ and $\ell$ are coprime.

{\em The case $\ell=3$.}
Let $\chi:\Ga_{\Q}\to\mu_4$ be a 1-cocycle
with the class $-3\in \Q^*/\Q^{*4}$. The restriction of $\chi$ to
$\Ga_{\Q(i)}$ is a homomorphism
$\chi':\Ga_{\Q(i)}\to \mu_4$. More precisely,
$\chi'$ is the quartic character $\Ga_{\Q(i)}\to\mu_4$
associated to the cyclic extension $\Q(i,\sqrt[4]{-3})$ of $\Q(i)$,
sending $\gamma\in \Ga_{\Q(i)}$ to
$\gamma(\sqrt[4]{-3})/\sqrt[4]{-3}$.

If $\pi$ is coprime to $3$, we can calculate $\chi'({\rm Frob}_\pi)$ 
using biquadratic residues
and the Gauss biquadratic reciprocity law (\cite[Thm. 4.21, Exercise 5.14]{Cox}):
$$\chi'({\rm Frob}_\pi)=\left(\frac{-3}{\pi}\right)_4=
\left(\frac{\pi}{3}\right)_4=\pi^2\ \in\ \mu_4 \subset \O/3.$$

The action of $z\in (\O/3)^*$ on $\End(E_3)^-$
sends an endomorphism $\varphi$ to $z\varphi z^{-1}=z\bar z^{-1}\varphi=z^{-2}\varphi$,
because for any $z\in (\O/3)^*=\F_9^*$ we have $\bar z=z^3$.
Thus $\Ga_{\Q(i)}$ acts on $\End(E_3)^-$ via the homomorphism
$\Ga_{\Q(i)}\to\mu_4$ that sends ${\rm Frob}_\pi$ to $\pi^{-2}$. 
We conclude that any isomorphism of $\mu_4$-modules 
$\End(E_3)^-\cong\O/3$ (which clearly exist) identifies the action of
$\Ga_{\Q(i)}$ on $\End(E_3)^-$ with the trivial action of $\Ga_{\Q(i)}$ on
$\O/3$ twisted by $\chi'^{-1}$.

Let $(\O/3)^{-1/3}$ be the twist of the $\Ga_\Q$-module $\O/3$ by $\chi^{-1}$.
We showed that there is a $\Ga_{\Q(i)}$-module isomorphism $f:\End(E_3)^-\to(\O/3)^{-1/3}$.
We claim that we can always compose $f$ with a $\mu_4$-module 
automorphism of $(\O/3)^{-1/3}$ to get a $\Ga_\Q$-module isomorphism.
The action of $\Ga_\Q$ on $(\O/3)^{-1/3}$ factors through the Galois group
$\Gal(\Q(i,\sqrt[4]{-3})/\Q)$, which is isomorphic to the dihedral 
group $D_8$. We obtain a homomorphism $\psi: D_8\to \GL(2,\F_3)$.
It restricts to an injection on $\Gal(\Q(i,\sqrt[4]{-3})/\Q(i))$, 
the normal cyclic subgroup of order 4, which implies that $\psi$ 
itself is injective.

We have the following for the subgroup structure of $\GL(2,\F_3)$:
Each cyclic subgroup, $C$, of order $4$ is contained in exactly one 
subgroup $H$ isomorphic to $D_8$.
Moreover any two elements of $H-C$ are conjugate in $\GL(2,\F_3)$ by 
an element which centralises $C$. These are enough to prove the claim.

The last sentence is clear from the construction.

\medskip

{\em The case $\ell=5$.}
Let $\chi:\Ga_{\Q}\to\mu_4$ be a 1-cocycle
with the class $5\in \Q^*/\Q^{*4}$. The restriction
$\chi':\Ga_{\Q(i)}\to \mu_4$ is the quartic character $\Ga_{\Q(i)}\to\mu_4$
associated to the cyclic extension $\Q(i,\sqrt[4]{5})$ of $\Q(i)$,
sending $\gamma\in \Ga_{\Q(i)}$ to
$\gamma(\sqrt[4]{5})/\sqrt[4]{5}$. We have $5=\theta\bar\theta$,
where $\theta=-1+2i$ and $\bar\theta=-1-2i$ are conjugate primary primes.
Then $\O/5\cong\O/\theta\oplus \O/\bar\theta$, where
$\O/\theta\cong \O/\bar\theta\cong\F_5$. The map $\O\to\O/\theta$
defines an isomorphism $\mu_4\to (\O/\theta)^*$, so an element
of $\mu_4\subset\O/5$ is uniquely determined by its projection to
$\O/\theta$.

Let $\pi$ be a primary prime in $\O$ which is coprime to 5.
Since $(N_{\Q(i)/\Q}(\theta)-1)/4=1$, the biquadratic reciprocity gives
$$\chi'({\rm Frob}_\pi)=\left(\frac{5}{\pi}\right)_4=
\left(\frac{\pi}{\theta}\right)_4\left(\frac{\pi}{\bar\theta}\right)_4=
\left(\frac{\pi}{\theta}\right)_4\left(\frac{\bar\pi}{\theta}\right)_4^{-1}=\pi\bar\pi^{-1}
\ \in \ \mu_4 \subset  (\O/5)^*.$$
By Gauss and Deuring ${\rm Frob}_\pi$ acts on $E_5$ as $\pi$, hence it acts
on $\End(E_5)^-$ by sending
an endomorphism $\varphi$ to $\pi\varphi \pi^{-1}=\pi\bar \pi^{-1}\varphi$.
We conclude that any isomorphism
of $\mu_4$-modules  $\End(E_5)^-\cong\O/5$ (which clearly exist)
identifies the action of
$\Ga_{\Q(i)}$ on $\End(E_5)^-$ with the trivial action of $\Ga_{\Q(i)}$ on
$\O/5$ twisted by $\chi'$.

One finishes the proof in the same way as for $\ell=3$, 
by looking at the subgroup structure of $\GL(2,\F_5)$.
$\Box$

\bpr \label{p3}
The $\Ga_\Q$-module $\Br(\ov{D})_3$ is isomorphic to the quartic twist
of the natural $\Ga_\Q$-module $\O/3$ by a cocycle with the class
$(-3abcd)^{-1} \in \Q^*/\Q^{*4}=\H^1(\Q,\mu_4)$.
In particular, if $-3abcd \notin \langle -4\rangle\subset\Q^*/\Q^{*4}$,
we have $\Br(\ov{D})_3^{\Ga_\Q}=\Br(\ov{D})_3^{\Ga_{\Q(i)}}=0$.
If $-3abcd \in \langle -4\rangle\subset\Q^*/\Q^{*4}$, then we have
$\Br(\ov{D})_3^{\Ga_{\Q(i)}}\cong (\Z/3)^2$ and
$\Br(\ov{D})_3^{\Ga_\Q}\cong \Z/3$.
\epr
{\em Proof.} The first statement follows from Lemma \ref{lem3}
and Proposition \ref{p2}. Since the kernel of the natural map
$\Q^*/\Q^{*4}\to \Q(i)^*/\Q(i)^{*4}$ is the cyclic subgroup of
order 2 generated by $-4$, everything apart from the last formula
is an immediate consequence. When
$-3abcd \in \Q^{*4}$ the last formula is clear. Finally, we need to
calculate the $\Ga_\Q$-invariants in the quartic twist of $\O/3$ by $-4$.
The twisted action of $\Ga_\Q$ factors through $\Gal(\Q(i)/\Q)$,
and the generator of this group sends $z$ to $-i\bar z$. It follows
that the rank of the invariant subgroup is 1. $\Box$

\bpr \label{p5}
The $\Ga_\Q$-module $\Br(\ov{D})_5$ is isomorphic to the quartic twist
of the natural $\Ga_\Q$-module $\O/5$ by a cocycle with the class
$5(abcd)^{-1} \in \Q^*/\Q^{*4}=\H^1(\Q,\mu_4)$.
In particular, if $125abcd \notin \langle -4\rangle\subset\Q^*/\Q^{*4}$,
we have $\Br(\ov{D})_5^{\Ga_\Q}=\Br(\ov{D})_5^{\Ga_{\Q(i)}}=0$.
If $125abcd \in \langle -4\rangle\subset\Q^*/\Q^{*4}$, then we have
$\Br(\ov{D})_5^{\Ga_{\Q(i)}}\cong (\Z/5)^2$ and
$\Br(\ov{D})_5^{\Ga_\Q}\cong \Z/5$.
\epr
{\em Proof.} This is proved in the same way as Proposition \ref{p3}.
$\Box$

\medskip

\noindent{\em End of proof of Theorem \ref{main}}. The statement of the
theorem is a formal consequence of  Lemma \ref{pp} and
Propositions \ref{c3}, \ref{p3} and \ref{p5}. $\Box$

\section{Proof of Theorem \ref{main2}} \label{section5}

\subsection{Preliminaries} \label{calc}

Without loss of generality we assume that $a,b,c,d$
are fourth power free non-zero integers such that $gcd(a,b,c,d)=1$.
In Proposition \ref{p1} we constructed a rational map of degree $32$ from
the diagonal quartic surface $D$ given by (\ref{*}) to the Kummer 
surface $Z=\Kum(E^{4ab}\times E^{4cd})$ with an affine equation
$$Y^2=(X^3-4abX)(T^3-4cdT).$$
If we denote $f=ax^4+by^4$, then
this map can be written as follows:
\begin{equation}
X=-4abx^2y^2f^{-1}, \ T=-4cdz^2w^2f^{-1}, \
Y=16abcdxyzw(ax^4-by^4)(cz^4-dw^4)f^{-3}. \label{map}
\end{equation}
Let $k$ be a field of characteristic zero, and let
$L=(x,y,z,w)\in D(k)$ such that $xyzwf\not=0$. Define $E'$
as the quadratic twist of $E^{4ab}$ by $f/x^2y^2$. Then
$E'$ is the quartic twist of $E$ by $4ab(f/x^2y^2)^2$,
that is
$$E': \quad u_1^2=t_1(t_1^2-4ab(f/x^2y^2)^2).$$
Then $P=(-4ab, 4abx^{-2}y^{-2}(ax^4-by^4))$ is a $k$-point in $E'$.
The quartic twist 
$$E'': \quad u_2^2=t_2(t_2^2-4cd(f/z^2w^2)^2)$$
is an elliptic curve over $k$
with a $k$-point $Q=(-4cd, 4cdz^{-2}w^{-2}(cz^4-dw^4))$.
Note that $E'\times E''$ is a quadratic twist of 
$E^{4ab}\times E^{4cd}$ defined over $k$, and hence there is a
natural rational map $\xymatrix{E'\times E'' \ar@{-->}[r]&Z\times_\Q k}$
of degree 2. This map sends $(P,Q)$ to the point $M\in Z(k)$ 
which is the image of $L\in D(k)$.
Explicitly this map is 
\begin{equation}
X=x^2y^2f^{-1}t_1, \quad T=z^2w^2f^{-1}t_2, 
\quad Y=(xyzw)^3f^{-3}u_1u_2. \label{map1}
\end{equation}

By Theorem \ref{main} we can assume that $\ell$ is 3 or 5. 
This assumption will be in force until the end of this paper.

Theorem \ref{main} implies, by (\ref{e2}), that if 
$\Br(\ov D)_\ell^{\Ga_\Q}\not=0$, then  
there is a non-zero homomorphism
of $\Ga_\Q$-modules $\varphi:E^{4cd}_\ell\to E^{4ab}_\ell$ such that
$[i]\varphi=-\varphi[i]$, which is unique up to multiplication
by an element of $\F_\ell^*$. It is easy to see that $\varphi$ is an
isomorphism. Indeed, $\Ker(\varphi)$ is $\mu_4$-invariant. For
$\ell=3$ the $\mu_4$-module $E^{4ab}_\ell\simeq\O/\ell$  is irreducible,
so $\varphi$ is an isomorphism in this case. For $\ell=5$ this
$\mu_4$-module
is the direct sum $\Ker[1+2i]\oplus\Ker[1-2i]$, however, neither factor
is a $\Ga_\Q$-submodule of $E^{4ab}_\ell$.

We need to recall from \cite[Section 3]{SZ2} 
how an element of $\Br(Z)_\ell$ 
is constructed from $\varphi$.
The multiplication by $\ell$ turns $E^{4ab}$ into an $E^{4ab}$-torsor
with structure group $E^{4ab}_\ell$. Let us call this torsor
$\T_1$, and similarly for the $E^{4cd}$-torsor $\T_2$.
The class $[\T_1]$ is an element of the \'etale cohomology
group $\H^1(E^{4ab},E^{4ab}_\ell)$. Similarly, 
$[\T_2]\in \H^1(E^{4cd},E^{4cd}_\ell)$. 
The homomorphism $\varphi$ gives rise to the
$E^{4cd}$-torsor $\varphi_*\T_2$ with structure 
group $E^{4ab}_\ell$. The Galois-equivariant Weil pairing
$$E^{4ab}_\ell\times E^{4ab}_\ell \lra \mu_\ell$$
gives rise to the cup-product pairing of \'etale cohomology groups
$$\H^1(E^{4ab}\times E^{4cd},E^{4ab}_\ell)\times 
\H^1(E^{4ab}\times E^{4cd},E^{4ab}_\ell) \lra 
\Br(E^{4ab}\times E^{4cd})_\ell.$$
Let us denote by $\sC\in \Br(E^{4ab}\times E^{4cd})_\ell$ 
the cup-product of the pullback of $\T_1$ 
via the projection to $E^{4ab}$ and the pullback 
of $\varphi_*\T_2$ via the projection to $E^{4cd}$.
According to \cite[Lemma 3.1, Prop. 3.3]{SZ2} the natural map
$$\Br(E^{4ab}\times E^{4cd})_\ell \to \Br(\ov E^{4ab}\times \ov E^{4cd})_\ell^{\Ga_\Q}=
\Hom_\Q(E^{4cd}_\ell,E^{4ab}_\ell)^-$$
sends $\sC$ to $\varphi$. By the proof of \cite[Thm. 2.4]{SZ2},
in particular, by diagram (16), the natural map $\Br(Z)_\ell\to
\Br(E^{4ab}\times E^{4cd})_\ell$ identifies $\Br(Z)_\ell$ with
the subgroup of $\Br(E^{4ab}\times E^{4cd})_\ell$ consisting
of the elements fixed by the antipodal
involution $[-1]$ on $E^{4ab}\times E^{4cd}$.
Let $\sB\in \Br(Z)_\ell$ be the element corresponding to $\sC$.

The quadratic twists of $E^{4ab}$ and $E^{4cd}$ 
by the same element $e\in k^*$ are elliptic curves
$E^{4abe^2}$ and $E^{4cde^2}$ over $k$. We then
obtain a homomorphism of $\Ga_k$-modules $E^{4cde^2}_\ell\to E^{4abe^2}_\ell$.
It is easy to see that the resulting element of $\Br(Z\times_\Q k)$
is precisely the image of $\mathcal B$ in $\Br(Z\times_\Q k)$.

Let $\sA\in \Br(D)_\ell$ be the image of $\sB$ under
the natural map $\Br(Z)_\ell\to\Br(D)_\ell$.

\bpr \label{coucou}
If $\nu$ is the automorphism of $D$ that alters the sign of
one of the coordinates and does not change the other three, then
$\nu^*\sA=-\sA$. Up to multiplication by an element
of $\F_\ell^*$, any permutation of variables 
in the equation of $D$ gives rise to the same element
$\sA\in \Br(D)_\ell$.
\epr
{\em Proof.} Without loss of generality we can assume
that $\nu(x,y,z,w)=(-x,y,z,w)$.
Let $[-1,1]$ be the automorphism of $E^{4ab}\times E^{4cd}$
that acts by $-1$ on $E^{4ab}$ and as the identity on $E^{4cd}$.
Since $[-1]^*[\T_1]=-[\T_1]$ we see that $[-1,1]^*\sC=-\sC$.
The involution $[-1,1]$ on $E^{4ab}\times E^{4cd}$ is compatible
with the involution $\kappa:Z\to Z$ that alters the sign of 
$Y$ and does not change $X$ and $T$, hence $\kappa^*\sB=-\sB$.
Finally, (\ref{map}) shows that the involutions $\nu$ and $\kappa$ 
are compatible as well, and therefore $\nu^*\sA=-\sA$. Any element of $\Br(D)_\ell$
that we construct after permuting the variables maps to a
generator of $\Br(\ov D)^{\Ga_\Q}_\ell\simeq\F_\ell$,
so, up to multiplication by an element of $\F_\ell^*$, the
difference of these elements is some $\alpha\in\Br_0(D)_\ell$.
We have $\alpha=\nu^*\alpha=-\alpha$, hence $\alpha=0$.
$\Box$

\bco \label{minus}
If $\ev_{\sA,p}$ is a constant map, it is identically zero.
\eco
{\em Proof.} By Proposition \ref{coucou} we have
$\sA(L)=\sA(\nu(L))=-\sA(L)\in\Br(\Q_p)_\ell$ for any $L\in D(\Q_p)$. 
But $\ell$ is odd, so $\sA(L)=0$. $\Box$

\medskip

\noindent{\bf Remark} In a similar way, using (\ref{map}) 
one checks that when $\mu_4$ is contained in $k$, 
the image of $(ix,y,z,w)\in D(k)$ in $Z(k)$ coincides with
the image of $([i]P,Q)$.

\medskip

Let $\chi$ be the homomorphism $E(k)\to \H^1(k,E_\ell)$, and
apply the same notation to the twists of $E$. Thus
the point $P\in E'(k)$ defines a $k$-torsor $\chi_P$ of $E'_\ell$,
and similarly $Q$ defines a $k$-torsor $\chi_Q$ of $E''_\ell$.
By the construction of $\sA$, $\sB$, $\sC$ we have
\begin{equation}
\sA(L)=\sB(M)=\chi_P\cup\varphi_*(\chi_Q) \in \Br(k)_\ell, \label{pairing}
\end{equation}
where
\begin{equation}\cup \ : \quad \H^1(k,E'_\ell) \ \times \ 
\H^1(k,E'_\ell) \ \lra
\ \Br(k)_\ell \label{cup}
\end{equation}
is the non-degenerate pairing induced by the Weil pairing
$E'_\ell\times E'_\ell \to \mu_\ell$. 

\medskip

Write $L\in D(\Q_p)$ as $(x,y,z,w)\in(\Z_p)^4$, where
$p$ does not divide all four coordinates. By the
implicit function theorem, up to replacing $L$ 
by its small deformation in $D(\Q_p)$ 
we can assume that $xyzw\not=0$ and $f(L)\not=0$.
Then $L$ goes to a well defined point $M\in Z(\Q_p)$, and 
$M$ lifts to the $\Q_p$-point $(P,Q)$ in $E'\times E''$,
where $E'$ and $E''$ are elliptic curves over $\Q_p$.

Given four non-zero elements $a,b,x,y$ of $\Z_p$ 
such that $f\not=0$ we denote the statement ``$P\in \ell E'(\Q_p)$" by 
$S_p(a,b:x,y)$.
Explicitly, it says that the point with $X$-coordinate $-4ab$
is divisible by $\ell$ on the elliptic curve 
$Y^2=X^3-4ab(f/x^2y^2)^2 X$. We shall use this statement in
the following equivalent form: 

\medskip

\noindent$S_p(a,b:x,y)\Leftrightarrow$ the point 
$(-4abx^2y^2,4abxy(ax^4-by^4))$ is divisible by $\ell$
on the elliptic curve $Y^2=X^3-4abf^2 X$ over $\Q_p$.

\medskip


By (\ref{pairing}) either $S_p(a,b:x,y)$ or $S_p(c,d:z,w)$ implies $\sA(L_p)=0$.

\subsection{Proof of (i)}

Let $m\in \Z_p$, $m\not=0$, and let $E^m$ be the elliptic curve 
$y^2=x^3-mx$ over $\Q_p$.
This Weierstrass equation is minimal if $\val_p(m)\leq 3$, see \cite[Remark
VII.1.1]{S}. Let $\tilde E^m$ be the (possibly, singular) curve given by
the equation $y^2=x^3-mx$ reduced modulo $p$.
To study the divisibility in $E^m(\Q_p)$ we recall the well known 
filtration
$$\ldots\subset E^m(\Q_p)_2\subset E^m(\Q_p)_1\subset E^m(\Q_p)_0 
\subset E^m(\Q_p),$$
where $E^m(\Q_p)_0$ is the subgroup of points that reduce to 
a point in $\tilde E^m_{\rm smooth}$, and $E^m(\Q_p)_1$ 
is the kernel of the reduction map $E^m(\Q_p)_0\to \tilde E^m_{\rm smooth}(\F_p)$.
One can also define $E^m(\Q_p)_n$ for $n\geq 1$ as the subgroup of $E^m(\Q_p)$
consisting of $0$ and the points $(x,y)$ such that $\val_p(x)\leq -2n$.

An application of Tate's algorithm \cite[IV.9.4]{S2} 
shows that at any prime $p$
the curve $E^m$ can only have good or additive reduction, and the order
of the quotient $E^m(\Q_p)/E^m(\Q_p)_0$ can only be 1, 2 or 4.
Thus for any prime $p$ we have an isomorphism
\begin{equation}
E^m_0(\Q_p)/\ell=E^m(\Q_p)/\ell. \label{neron}
\end{equation}
The quotient $E^m(\Q_p)_0/E^m(\Q_p)_1$ is naturally isomorphic to 
$\tilde E^m_{\rm smooth}(\F_p)$ \cite[Prop. VII.2.1]{S}. 
Pannekoek \cite[Thm. 1, Lemma 9]{P}, using the well known
theory of the formal group attached to an elliptic curve 
\cite[Ch. IV]{S}, proves that in the case of additive reduction
$E^m(\Q_p)_0$ is topologically isomorphic to $\Z_p$ or to $p\Z_p\times\F_p$,
where the last case occurs if and only if $p=5$ and $m\equiv 15\bmod 25$.

\ble \label{two}
The group $E^m(\Q_2)$ is divisible by any odd prime.
Thus $S_2(a,b:x,y)$ holds for any non-zero $a,b,x,y\in\Z_2$.
\ele
{\em Proof.} Tate's algorithm shows that 
$E^m$ has additive reduction for any $m\in\Q_2^*$.
The lemma follows from (\ref{neron}) and the above remarks. $\Box$

\medskip

From now on assume $p>2$. Then the theory of
the formal group attached to $E^m$ gives compatible isomorphisms
$E^m(\Q_p)_n\simeq p^n\Z_p$ for $n\geq 1$, see \cite[Thm. IV.6.4(b)]{S}.

\ble\label{lc5fo}
Let $p\not=\ell$ be an odd prime. If $a,b,x,y$ are in $\Z_p^*$ and 
$\val_p(f)\geq 1$, then  $S_p(a,b:x,y)$ holds.
\ele
{\em Proof.} If $\val_p(f)$ is odd, then $E'$ has additive reduction at $p$.
By (\ref{neron}) the group $E'(\Q_p)$ is divisible by $\ell$.

Suppose $f=up^{2n}$, where $u\in\Z_p^*$ and $n\geq 1$. The coordinate change
$X'=p^{-2n}X$, $Y'=p^{-3n}Y$ reduces the equation of $E'$ to an
equation the reduction modulo $p$ of which defines a smooth curve 
$\tilde E'$ over $\F_p$. The $X'$-coordinate of $P$ has valuation $-2n$,
and so $P$ is in $E'(\Q_p)_1$. This is a pro-$p$-group, 
hence $P\in \ell E'(\Q_p)$. $\Box$

\bpr\label{Pr5fo}
Let $D=[a,b,c,d]$, where $a,b,c,d\in\Q^*$. If $p\not=\ell$, then
the map $\ev_{\sA,p}$ is zero.
\epr
{\em Proof.} By Corollary \ref{minus} it is enough to show that
the map $\ev_{\sA,p}$ is constant. If $p$ does not divide $2abcd$,
this was proved in \cite[Cor. 3.3]{CS1}.

Since $S_2(a,b:x,y)$ holds by Lemma \ref{two}, we have
$\sA(L)=0$ for any $L\in D(\Q_2)$.

Now let $p$ be an odd prime. Since $p\not=\ell$, we see from Theorem
\ref{main} that $4$ divides $\val_p(a)+\val_p(b)+\val_p(c)+\val_p(d)$.
In the proof that $\ev_{\sA,p}$ 
is a constant map we are allowed to rename
the variables in the equation of $D$. By doing so
and by multiplying $a,b,c,d$ by elements of $\Q^{*4}$
we can assume that $r=[\val_p(a),\val_p(b),\val_p(c),\val_p(d)]$
is either $[0,0,2,2]$ or $[0,0,1,3]$.
If $r=[0,0,2,2]$, then either both $x$ and $y$ are units,
or both $z$ and $w$ are units.
By Lemma \ref{lc5fo} either $S_p(a,b:x,y)$ or 
$S_p(c/p^2,d/p^2:z,w)$ holds, and the last property is
equivalent to $S_p(c,d:z,w)$.
If $r=[0,0,1,3]$, then $S_p(a,b:x,y)$ holds by the same lemma. 
$\Box$

\medskip

Theorem \ref{main2} (i) is now proved, because 
$\sA$ generates $\Br(D)_\ell$ modulo $\Br_0(D)_\ell$. 

The following corollary follows from Proposition \ref{Pr5fo} and
the reciprocity law; it will be used in the next section.

\bco \label{global}
Let $D=[a,b,c,d]$, where $a,b,c,d\in\Q^*$. Then for any $P\in D(\Q)$
we have $\ev_{\sA,\ell}(P)=0$.
\eco

\subsection{Proof of (ii) for $\ell=5$}

\subsubsection{Computation of the pairing}

In this section we write $\Ga=\Ga_{\Q_5}$ and denote by
$E^m$ the elliptic curve $y^2=x^3-mx$ where $m\in \Q_5^*$.

Since $-1\in\Z_5^{*2}$ we can define $i\in \Z_5$
by the property that $1+2i$ is a generator of the maximal ideal
$(5)\subset\Z_5$, or, equivalently, $i\equiv 2\bmod 5$.
Then $1-2i\in\Z_5^*$.
We have a natural decomposition of $\Ga$-modules
$E^m_5=E^m_{1+2i} \oplus E^m_{1-2i}$, and
the induced decomposition
$$\H^1(\Q_5,E^m_5)=\H^1(\Q_5,E^m_{1+ 2i})\oplus \H^1(\Q_5,E^m_{1-2i}).$$
Since the restriction of the skew-symmetric Weil pairing to
$E^m_{1\pm 2i}$ is trivial, each of the subspaces
$\H^1(\Q_5,E^m_{1\pm 2i})$ is isotropic.
By the non-degeneracy of the $\cup$-product, these subspaces
are maximal isotropic subspaces of $\H^1(\Q_5,E^m_5)$,
each of dimension $\frac{1}{2}\dim \H^1(\Q_5,E^m_5)$.
We conclude that (\ref{cup}) induces a non-degenerate pairing
\begin{equation}
\H^1(\Q_5,E^m_{1+2i})\times \H^1(\Q_5,E^m_{1-2i})\lra \Br(\Q_5)_5
\cong\frac{1}{5}\Z/\Z. \label{cup1}
\end{equation}
Recall that $E^m(\Q_5)/5$ is a maximal isotropic subspace of
$\H^1(\Q_5,E^m_5)$. 
The Chinese remainder theorem gives an isomorphism
$$E^m(\Q_5)/5=E^m(\Q_5)/[1+2i] \oplus E^m(\Q_5)/[1-2i].$$
We have the following commutative square of inclusions:
$$\begin{CD}
  E^m(\Q_5)/[1+2i] @>>> \H^1(\Q_5,E^m_{1+2i}) \\
      @VV[1-2i] V             @VVV     \\
  E^m(\Q_5)/5 @>>> \H^1(\Q_5,E^m_{5})
\end{CD}$$
where the horizontal maps come from the corresponding Kummer sequences
and the right vertical map is induced by the inclusion
$E^m_{1+2i}\to E^m_5$.

\bpr \label{pr-pr}
If $m\in\Q_5^*$ is not in $(1\pm 2i)\Q_5^{*4}$, then
$E^m(\Q_5)/[1-2i]=0$ and 
$$E^m(\Q_5)/5=E^m(\Q_5)/[1+2i]=\H^1(\Q_5,E^m_{1+2i})$$
is a $1$-dimensional $\F_5$-vector subspace of 
the $2$-dimensional space $\H^1(\Q_5,E^m_{5})$.
\epr
{\em Proof.} Let $\mathcal F$ be the formal group defined by
$E^m$, see \cite[Example IV.2.2.3]{S}. One can consider the abelian group
$\mathcal F(5\Z_5)$ on the set $5\Z_5$ with the group law
defined by $\mathcal F$, see \cite[Example IV.3.1.3]{S}.
By \cite[VII.2.2]{S} the map sending $(x,y)$
to $-x/y$ is an isomorphism $E^m(\Q_5)_1\to \mathcal F(5\Z_5)$.
This isomorphism is clearly compatible with the action of $[i]$
on $E^m(\Q_5)_1$ and the action of $i\in \Z_5^*$ on the set $5\Z_5$.
Since the residual characteristic is not 2 we can apply
\cite[Thm. IV.6.4 (b)]{S} which says that the formal logarithm
defines an isomorphism of abelian groups $\mathcal F(5^n\Z_5)\to
5^n\Z_5$ for any $n\geq 1$. The composed isomorphism $E^m(\Q_5)_1\to 5\Z_5$
translates the action of $[i]$
on $E^m(\Q_5)_1$ into the action of $i\in \Z_5^*$ on $5\Z_5$.

Write $m=5^n u$, where $n\in\Z$, $0\leq n\leq3$, and $u\in\Z_5^*$.
First suppose that $m\in \Z_5^*$, then $E^m$ has good reduction.
Our assumption $m\notin(1-2i)\Z_5^{*4}$ implies that $m$
is not congruent to $2$ modulo $5$.
Let us denote by $\tilde E^m$ the elliptic curve over $\F_5$
which is the reduction of $E^m$. An elementary calculation
shows that for such values of $m$, the $5$-torsion of 
$\tilde E^m(\F_5)$ is trivial, hence $[1\pm2i]$-torsion is trivial too.
Since $1-2i\in\Z_5^*$, the endomorphism $[1-2i]$ is invertible on
$E^m(\Q_5)$. 

Now let $m\in 5^n\Z_5^*$, where $n=1,2,3$. Then $E^m$ has additive reduction.
R.~Pannekoek shows in
\cite[Def. 10]{P} that in this situation the formal group $\mathcal F$
defines an abelian group $\mathcal F(\Z_5)$ on the set $\Z_5$.
He then proves \cite[Prop. 11 (1)]{P} that the map sending $(x,y)$
to $-x/y$ is an isomorphism $E^m(\Q_5)_0\to \mathcal F(\Z_5)$.
As before, this isomorphism is obviously compatible with the action of $[i]$
on $E^m(\Q_5)_0$ and the action of $i\in \Z_5^*$ on $\mathcal F(\Z_5)$.
By \cite[Prop. 17]{P} there is a natural isomorphism
between $\mathcal F(\Z_5)$ and $\Z_5$ with its usual group structure, unless
$m\equiv 15\equiv 1+2i \bmod 25$. This congruence easily implies $m\in (1+2i)\Z_5^{*4}$
which is excluded by our assumption.
The number of connected components of the N\'eron model
of an elliptic curve with additive reduction is at most 4.
So the endomorphism $[1-2i]$ is invertible on
$E^m(\Q_5)$ in this case as well. 

Thus in all cases we have $E^m(\Q_5)/[1-2i]=0$ and
$E^m(\Q_5)/5=E^m(\Q_5)/[1+2i]$ is a 1-dimensional vector space over $\F_5$.
Since $E^m(\Q_5)/[1+2i]\subset \H^1(\Q_5,E^m_{1+2i})$ and the
dimension of the last space equals the dimension of $E^m(\Q_5)/5$,
we see that the injective image of $E^m(\Q_5)/5$
in $\H^1(\Q_5,E^m_5)$ is $\H^1(\Q_5,E^m_{1+2i})$. $\Box$

\medskip

If $\varphi:E^{m_2}_5\to E^{m_1}_5$ is an isomorphism of
$\Ga$-modules which anti-commutes with $i$, then
$\varphi=\varphi'+\varphi''$, where $\varphi'$ is an isomorphism
of $\Ga$-modules $E^{m_2}_{1+2i}\to E^{m_1}_{1-2i}$,
and $\varphi''$ is an isomorphism
of $\Ga$-modules $E^{m_2}_{1-2i}\to E^{m_1}_{1+2i}$.
Therefore, the induced isomorphism
$\varphi_*:\H^1(\Q_5,E^{m_2}_5)\tilde\lra \H^1(\Q_5,E^{m_1}_5)$
is the sum $\varphi'_*+\varphi''_*$, where
$$\varphi'_*:\H^1(\Q_5,E^{m_2}_{1+2i})\tilde\lra \H^1(\Q_5,E^{m_1}_{1-2i}),$$
and similarly for $\varphi''_*$. Now the non-degeneracy of (\ref{cup1})
implies the non-degeneracy of the pairing
\begin{equation}
\H^1(\Q_5,E^{m_1}_{1+2i})\times \H^1(\Q_5,E^{m_2}_{1+2i})\lra \Br(\Q_5)_5
\cong\frac{1}{5}\Z/\Z \label{cup2}
\end{equation}
given by $x\cup \varphi'_*(y)$.

\bco\label{calc5}
Let $m_1,m_2\in\Q_5^*$ be such that $m_1m_2\in 5\Q_5^{*4}$, where
$m_1$ and $m_2$ are not in $(1\pm 2i)\Q_5^{*4}$.
Let $\varphi:E^{m_2}_5\to E^{m_1}_5$ be an isomorphism of
$\Ga_{\Q_5}$-modules that anti-commutes with $[i]$.
For any $P\in E^{m_1}(\Q_5)$ and $Q\in E^{m_2}(\Q_5)$ we have

{\rm (1)} \quad $\chi_P\cup\varphi_*(\chi_Q)=0$ if and only if
$P\in 5E^{m_1}(\Q_5)$ or $Q\in 5E^{m_2}(\Q_5)$;

{\rm (2)} \quad $\chi_{[i]P}\cup\varphi_*(\chi_Q)=\chi_P\cup\varphi_*(\chi_{[i]Q})=
2\,\chi_P\cup\varphi_*(\chi_Q).$

\eco
{\em Proof.} (1) Applying Proposition \ref{pr-pr} to $m_1$
and $m_2$, we obtain 
$\chi_P\cup\varphi_*(\chi_Q)=x\cup \varphi'_*(y)$
for some non-zero elements $x\in\H^1(\Q_5,E^{m_1}_{1+2i})$ and 
$y\in\H^1(\Q_5,E^{m_2}_{1+2i})$. The last pairing is 
a non-degenerate pairing of 1-dimensional spaces, hence
we prove (1).

(2) The automorphism $[i]$ acts on $E^m_{1+2i}$, and hence also
on $\H^1(\Q_5,E^m_{1+2i})$, as the multiplication by $2$.
We conclude by Proposition \ref{pr-pr}. $\Box$

\bco\label{lc5f5}
Let $D=[a,b,c,d]$, where $a,b,c,d\in\Q_5^*$ are subject to
conditions $5^3abcd\in\Q_5^{*4}$ and $ab\in \Q_5^{*2}$. If there exists
a point $(x,y,z,w)\in D(\Q_5)$ such that neither $S_5(a,b:x,y)$
nor $S_5(c,d:z,w)$ holds, then $\ev_{\sA,5}$ is not the zero map.
\eco
{\em Proof.} Since $ab\in \Q_5^{*2}$, we have 
$4abf^2\notin (1\pm2i)\Q_5^{*2}$, hence
$4cdf^2\notin (1\pm 2i)\Q_5^{*4}$. Now
Corollary \ref{calc5} (1) implies the lemma. $\Box$

\subsubsection{Diagonal quartic surfaces over $\Q_5$}

Consider the following diagonal quartic surfaces defined over $\Q$:
$$A_n=[1,-1,n,-5n^3], \quad B_n=[1,-1,5^2n,-5^3n^3], \quad C=[2,2,4,5],$$
where $n$ is an integer not divisible by $5$. Each of these surfaces 
has an obvious $\Q$-point and satisfies the property 
$ab\equiv -1\bmod 5$, hence $ab\in \Q_5^{*2}$.
In the following lemma and thus everywhere in this section it is enough 
to consider $n=1$ or $2$, though we shall not use this.

\ble\label{ex5}
Let $a,b,c,d\in\Q_5^*$ be such that $abcd \in 5\Q_5^{*4}$ and 
$D(\Q_5)\neq \emptyset$.
By permuting the variables and multiplying the coefficients
$a,b,c,d$ by a common constant
and by elements of $\Q_5^{*4}$, the surface $D$ can be reduced
to one of the surfaces $A_n$, $B_n$ or $C$.
In particular, we can assume without loss of generality that $ab\in\Q_5^{*2}$.
\ele
{\em Proof.} Let $r=[\val_5(a),\val_5(b),\val_5(c),\val_5(d)]$. 
We can assume that $r$ is $[0,0,0,1]$ or $[0,0,2,3]$.
If $r=[0,0,0,1]$, then the congruence $ax^4+by^4-cz^4\equiv 0\bmod 5$
has a non-trivial solution $(x_0,y_0,z_0)$. If $5|x_0y_0z_0$,
say $5|z_0$, we reduce $D$ to $A_n$; otherwise
we reduce it to $C$. When $r=[0,0,2,3]$ we have 
$a+b\equiv 0 \mod 5$, and then $D$ reduces to $B_n$. $\Box$

\ble \label{goodred}
Assume that $ab\in \Z_5^{*2}$.
If $\val_5(f)=0$ and $r=\val_5(xy)>0$,
then $S_5(a,b:x,y)$ holds if and only if $r>1$.
\ele
{\em Proof.} In this case the curve $E'$ given by $Y^2=X^3-4abf^2 X$
has good reduction. The point $P$ with $X$-coordinate $-4abx^2y^2$ 
reduces to the 2-torsion point $(0,0)$. Since $(0,0)=5(0,0)$,
the point $P$ is divisible by $5$ in $E'(\Q_5)$
if and only if $P_1=P+(0,0)$ is.
The $X$-coordinate of $P_1$ is $f^2/x^2y^2$, and
$\val_5(X(P_1))=-2r$, so that $P_1$ is in $E'(\Q_5)_1\simeq 5\Z_5$.
The assumption $ab\in \Z_5^{*2}$ implies that the group
of points on the reduced curve has no 5-torsion. Thus
$P_1\in 5E'(\Q_5)$ if and only if $P_1\in 5E'(\Q_5)_1$.
Since $\val_5(X(P_1))=-2r$, we conclude 
that $P_1\in 5^r\Z_5$, $P_1\notin 5^{r+1}\Z_5$.
In other words, $P_1$ is divisible by $5^{r-1}$ 
but not by $5^r$ in $E'(\Q_5)$. $\Box$

\bpr\label{Pr5f5}
Let $D=[a,b,c,d]$, where $a,b,c,d\in\Q^*$ are such that
$5^3abcd$ is in the subgroup $\langle-4\rangle\Q^{*4}$
and $D(\Q_5)\neq \emptyset$. Then for any element in
$\Br(D)_5$ that is not in $\Br_0(D)$ the evaluation map
$D(\Q_5)\to \frac{1}{5}\Z/\Z$ is surjective.
\epr
{\em Proof.} By Lemma \ref{ex5} we may assume that $D$ is one of the 
surfaces $A_n$, $B_n$ or $C$, in particular,
$ab\in\Q_5^{*2}$. It is enough to prove the surjectivity 
of $\ev_{\sA,5}$ where $\sA\in\Br(D)$ is constructed in Section \ref{calc}.
Each of our surfaces is defined over $\Q$
and contains a $\Q$-point. By Corollary \ref{global} for
any $\Q$-point $P$ we have $\ev_{\sA,5}(P)=0$.
We shall exhibit a $\Q_5$-point $L=(x,y,z,w)$ 
on each of our surfaces such that
neither $S_5(a,b:x,y)$ nor $S_5(c,d:z,w)$ holds.
Then $\sA(L)\neq0$ by Corollary \ref{lc5f5}. 
Corollary \ref{calc5} (2) and the remark after Corollary \ref{minus}
imply that the values of $\ev_{\sA,5}$ at the $\Q_5$-points
$(i^nx,y,z,w)$ for $n=0,1,2,3$ give all non-zero elements of 
$\frac{1}{5}\Z/\Z$. This will prove the proposition. 

$A_n$. {\em In this case we prove more, namely, 
if $25$ does not divide $xyzw$,
then neither $S_5(a,b:x,y)$ nor $S_5(c,d:z,w)$ holds.}
An example of such a point is $(\alpha,1,5,1)$ for an appropriate $\alpha\in\Z_5^*$.

Indeed, suppose that $x,y,z,w\in\Z_5$, not all divisible
by $5$, are the coordinates of a point in $D(\Q_5)$ such that 
$\val_5(xyzw)\leq 1$. Then $\val_5(z)$ is 0 or 1.

Let us first assume that $\val_5(z)=1$. This implies $x,y,w\in\Z_5^*$,
and hence $\val_5(f)=1$. Then the curve $E'$ given by $Y^2=X^3-4abf^2 X$
has additive reduction. The point $P$ with $X$-coordinate 
$-4abx^2y^2\in\Z_5^*$ is in $E'(\Q_5)_0$ but not in $E'(\Q_5)_1$,
so it is not divisible by $5$ in $E'(\Q_5)$. Thus
$S_5(a,b:x,y)$ does not hold. The curve $E''$ given by $Y^2=X^3-4cdf^2 X$
also has additive reduction. The point $Q$ has $X$-coordinate 
$-4cdz^2w^2$ of valuation $3$, hence it reduces to the singular point
$(0,0)$ of the reduction. The $X$-coordinate of $Q+(0,0)$ is $f^2/z^2w^2\in\Z_5^*$,
so $Q+(0,0)$ is not divisible by $5$ in $E''(\Q_5)$. Since $5(0,0)=(0,0)$,
the property $S_5(c,d:z,w)$ does not hold.

Now let $\val_5(z)=0$. This rules out the possibility that $x,y\in\Z_5^*$,
so the valuation of $x$ or $y$ is 1, and hence $w\in\Z_5^*$.
We now have $\val_5(f)=0$.
Since $-4abf^2\in\Z_5^*$ the curve $E'$ has good reduction,
and by Lemma \ref{goodred} the property $S_5(a,b:x,y)$ does not hold.
The curve $E''$ given by $Y^2=X^3-4cdf^2 X$
has additive reduction. The point $Q$ with $X$-coordinate 
$-4cdz^2w^2$ reduces to $(0,0)$. The point $Q+(0,0)$ has
$X$-coordinate $f^2/z^2w^2\in\Z_5^*$, so $Q+(0,0)$, and hence
also $Q$, is not divisible by $5$ in $E''(\Q_5)$, so
the property $S_5(c,d:z,w)$ does not hold. 
This finishes the proof of our claim.

\medskip

$B_n$. {\em In this case we prove more, namely, 
if $5$ does not divide $zw$,
then neither $S_5(a,b:x,y)$ nor $S_5(c,d:z,w)$ holds.}
An example of such a point is $(\alpha,1,1,1)$ 
for an appropriate $\alpha\in\Z_5^*$.

Suppose $x,y,z,w\in\Z_5$ are the coordinates of a point 
in $D(\Q_5)$ such that $\val_5(z)=\val_5(w)=0$.
Then $\val_5(f)=2$. The curve $E'$ has good reduction,
but its Weierstrass equation is not minimal. Changing variables to
make the equation minimal we see that $P$ is in $E'(\Q_5)_1$
but not in $E'(\Q_5)_2=5E'(\Q_5)_1$. Since 
the group of points on the reduced curve has no 5-torsion
we conclude that $S_5(a,b:x,y)$ does not hold. The 
equation $Y^2=X^3-4cdf^2X$ of the curve 
$E''$ is not minimal. Write $X=5^4X'$, $Y=5^6Y'$,
then the equation becomes $Y'^2=X'^3-4cdf^2 5^{-8}X'$.
The reduction is additive, and the $X'$-coordinate of $Q$
is $-4cdz^2w^2 5^{-4}$, so $Q$ reduces to the singular point
$(0,0)$. The $X'$-coordinate of $Q+(0,0)$ is $5^{-4}f^2/z^2w^2$,
so this point is in $E''(\Q_5)_0$ but not in $E''(\Q_5)_1$.
Thus $Q+(0,0)$, and hence also $Q$, is not divisible by $5$
in $E''(\Q_5)$, so $S_5(c,d:z,w)$ does not hold. 

\medskip

$C$. Let $L=(1,2,z,1)$, where $z\in\Z_5^*$. The equation of $E'$ is
$Y^2=X^3-16\cdot (34)^2 X$, so $E'$ has good reduction.
We need to prove that $P=(-64,-960)$ is not divisible by $5$
in $E'(\Q_5)$.
By writing $X=4X'$, $Y=8Y'$ we reduce the equation to
$Y'^2=X'^3-(34)^2 X'$. In new coordinates $P$
is the point $(-16,-120)$.
The reduction of this point coincides with the reduction
of the 2-torsion point $M=(34,0)$. Using {\tt sage}
one sees immediately that the valuation of the $X$-coordinate
of $P+M$ is $-2$. Hence $P-M\in E'(\Q_5)_1$, but $P-M$
is not in $5E'(\Q_5)_1$. The number of
$\F_5$-points on the reduction of $E'$ is 8, hence is prime to 5.
It follows that $P\in 5E'(\Q_5)$ if and only if 
$P-M\in 5E'(\Q_5)_1$. Thus $P$ is not divisible by 5 in $E'(\Q_5)$,
so $S_5(2,2:1,2)$ does not hold.

The curve $E''$ with equation $Y^2=X^3-4cdf^2X$ has additive reduction.
The valuation of the $X$-coordinate of the point $Q+(0,0)$ is 
$\val_5(f^2/z^2w^2)=0$, so this point is not in 
$E''(\Q_5)_1$. Hence $Q$ is not divisible by $5$, so that 
$S_5(4,5:z,1)$ does not hold. We conclude that $\sA(L)\not=0$.
$\Box$

\medskip

This finishes the proof of Theorem \ref{main2} (ii) for $\ell=5$.

\subsection{Proof of (ii) for $\ell=3$}

\subsubsection{Computation of the pairing}

Let $E^m$ be the elliptic curve $y^2=x^3-mx$, where $m\in \Q_3^*$.

\bpr\label{Bpr3}
Suppose that $-3m_1m_2 \in \langle-4\rangle\Q_3^{*4}$ 
and let $\varphi:E^{m_2}_3\to E^{m_1}_3$ be an isomorphism of $\Ga_{\Q_3}$-modules.
For any $P\in E^{m_1}(\Q_3)$ and $T\in E^{m_2}(\Q_3)$ we have
$\chi_P\cup\varphi_*(\chi_T)=0$ if and only if 
$P\in 3E^{m_1}(\Q_3)$ or $T\in 3E^{m_2}(\Q_3)$.
\epr
{\em Proof.} One implication being trivial, we suppose that 
$\chi_P\neq 0$ and $\chi_T\neq 0$.
The image of $E(\Q_3)/3$ under the Kummer map is a  
1-dimensional maximal isotropic subspace of $\H^1(\Q_3,E_3)$. 
Therefore it is enough to show that $\chi_P\neq \pm \varphi(\chi_T)$
in $\H^1(\Q_3,E^{m_1}_3)$. 

The elliptic curves $E=E^{m_1}$ and $C=E^{m_2}$ are isomorphic over 
$F=\Q_3(i,\sqrt[4]{3})$, but not over $\Q_3(i,\sqrt{3})$. Let 
$\lambda:C\to E$ be an isomorphism over $F$,
and let $Q=\lambda(T)$. Write $L=\Q_3(E_3)=\Q_3(C_3)$. 
Then $L/F$ is a ramified quadratic extension and
$L/\Q_3$ is a Galois extension.

We now prove a series of lemmas.

\ble\label{li}
If $P\not\in 3E^{m_1}(\Q_3)$ and $T\not\in 3E^{m_2}(\Q_3)$, then 
the images of $P,[i]P,Q,[i]Q$ in $E(L)/3$ are linearly independent.
\ele
{\em Proof.} 
Since the restriction map $H^1(F,E_3)\to H^1(L,E_3) $ is injective, 
it suffices to show linear independence in $E(F)/3$.
The map $H^1(\Q_3,E_3)\to H^1(F,E_3) $ is also injective, 
so $P$ does not go to zero in $E(F)/3$.

If $[i]P= \pm P$ in $E(F)/3$, then $2P=0$ 
which is impossible. Therefore $P$ and $[i]P$ are linearly 
independent. The same argument shows that $Q$ and $[i]Q$
are linearly independent in $E(F)/3$.

Note that $\lambda$ identifies $C(F)$ with $E(F)$ twisted by $-1$ 
as $\Gal(F/\Q_3(i,\sqrt{3}))$-modules.
Therefore the non-trivial element of $\Gal(F/\Q_3(i,\sqrt{3}))$ acts by
$1$ on the span of $\{P,[i]P\}$ and by $-1$ on the span of 
$\{Q,[i]Q\}$. The result follows. $\Box$

\medskip

For a field extension $k$ of $\Q_3$ and 
a point $R=(x_R,y_R)\in E^m(k)$ such that $6R\neq0$, consider
the polynomial
$$
g_R(t)=t^9+12mt^7+30m^2t^5-36m^3t^3+9m^4t-x_R(3t^4-6mt^2-m^2)^2.
$$
Let $k_R$ be the splitting field of $g_R(t)$.
The roots of $g_R(t)$ are the $x$-coordinates of the points 
$R'$ such that $3R'=R$. Note that $R'$ is defined over $k(x_{R'})$,
because otherwise $2R=0$.

\ble\label{lex}
Let $k$ be a field extension of $\Q_3$. Let $R\in E(k)$, $6R\neq0$.

{\rm (1)} If an extension $k$ contains two roots of $g_R(t)$, 
then $L\subset k(i)$.

{\rm (2)} Suppose $k\subset L$. Then $L_R/k$ is a Galois extension.
For any root $r$ of $g_R(t)$ we have $L_R=L(r)$. 
Moreover, $[L_R:L]$ divides $9$.

{\rm (3)} If $[L_R:L]=3$, then $\Gal(L_R/F)\cong D_6$.
\ele
{\em Proof.} (1) If two roots of $g_R(t)$ are in $k$, then 
a non-zero element of $E_3$ is defined over $k$ and so $E_3\subset E(k(i))$.

(2) The extension $L_R/k$ is Galois because it is the composite of 
the Galois extensions $L/k$ and $k_R/k$. 
Since $E_3\subset E(L)$, we also have the second assertion. 
The last assertion follows because the irreducible 
factors of $g_R(t)$ in $L[t]$ must have the same degree.

(3) By part (2) we see that $L_R/F$ is a Galois extension of degree 6. 
Therefore it suffices to show that it is not abelian.
It is easy to see that at least one irreducible factor of $g_R(t)$ in 
$F[t]$ must have degree $3$, call it $p(t)$.
If $L_R/F$ were abelian, then every subextension would be Galois. 
Hence by adjoining a root of $p(t)$ we would obtain an extension 
of $F$ of degree 3 containing two roots of $g_R(t)$. 
By part (1) this extension would contain $L$, which 
is a contradiction. $\Box$

\ble\label{uneq}
We have $L_Q\neq L_P$.
\ele
{\em Proof.} Assume $L_Q=L_P$. The $\F_3$-vector space
$$
\H^1(\Gal(L_P/L),E_3)=\Hom(\Gal(L_P/L),E_3)=\Ker[\H^1(L,E_3)\to \H^1(L_P,E_3)]
$$
has dimension at least $4$ by Lemma \ref{li}. Since $|\Gal(L_P/L)|$ divides
$9$ by Lemma \ref{lex}, we deduce that $\Gal(L_P/L)\cong \Z/3 \times \Z/3$.

Hence the kernel of the map $\H^1(L,E_3)\to \H^1(L_P,E_3)$ 
is the 4-dimensional linear span of 
$P,[i]P,Q,[i]Q$ in $E(L)/3$.
Let $k$ be any subextension of $L_P/L$ of degree 3.
The kernel of $\H^1(L,E_3)\to \H^1(k,E_3)$ is the 2-dimensional span of 
$R,[i]R$ for some $R$ in the span of $P,[i]P,Q,[i]Q$. Thus $k=L_R$.
By Lemma \ref{lex}, $k$ is a Galois 
extension of $F$ with the Galois group isomorphic to $D_6$.
But then $\Gal(L_P/F)$ contradicts the following 
fact from group theory: there does not exist a group 
of order $18$ with at least $4$ normal subgroups of order $3$ 
each of which has a quotient isomorphic to $D_6$.
The easy proof of this fact is left to the reader. $\Box$

\medskip

\noindent{\em End of proof of Proposition \ref{Bpr3}.}
Let us choose $P'$ such that $3P'=P$ and $Q'$ such that $3Q'=Q$. 
Since $L_Q \neq L_P$ by Lemma \ref{uneq}, we can find an element
$g \in \Gamma_L$ that fixes exactly one of $P'$ and $Q'$. 
Therefore exactly one of $\varphi(\chi_T)(g)$ and $\chi_P(g)$ 
is non-zero. Hence $\varphi(\chi_T)\neq \pm \chi_P$ in $\H^1(L,E_3)=\Hom(\Gamma_L,E_3)$.
The proposition is proved. $\Box$

\medskip

In the following immediate corollary
$\sA$ is the element of $\Br(D)$ constructed in Section \ref{calc}.

\bco \label{co2}
Let $D=[a,b,c,d]$, where $a,b,c,d\in\Q^*$ are such that
$-3abcd$ is in $\langle-4\rangle\Q^{*4}$.
If there is a point $(x,y,z,w)\in D(\Q_3)$ 
such that neither $S_3(a,b:x,y)$ nor 
$S_3(c,d:z,w)$ holds, then $\ev_{\sA,3}$ is a non-zero map.
\eco

\subsubsection{Diagonal quartic surfaces over $\Q_3$}

\ble\label{ex3}
Any diagonal quartic surface $D=[a,b,c,d]$ over $\Q_3$ such that
$-3abcd$ is in $\langle-4\rangle\Q_3^{*4}$ 
and $D(\Q_3)\neq \emptyset$,
can be reduced by a permutation of variables and by multiplication of 
the coefficients $a,b,c,d$ by elements of $\Q_3^{*4}$ 
and by a common constant in $\Q_3^*$, to one of the following surfaces:
$$A=[1,1,1,\pm 27], \quad B=[1,-1,3,\pm 9], \quad C=[1,1,2,\pm 27].$$
\ele
{\em Proof}. We can assume without loss of generality that
$[val_3(a),val_3(b),val_3(c),val_3(d)]$ equals $[0,0,0,3]$ or $[0,0,1,2]$.
In the first case if $a,b,-c$ are equal modulo 
$3$, then we reduce to $C$; otherwise our surface reduces to $A$.
In the second case the solubility in $\Q_3$ 
implies $a+b\equiv 0 \mod 3$ and we can take $a=1$.
This gives $B$.  $\Box$

\bpr\label{Pr3f3}
Let $D=[a,b,c,d]$, where $a,b,c,d\in\Q^*$ are such that
$-3abcd$ is in $\langle-4\rangle\Q^{*4}$ and $D(\Q_3)\neq \emptyset$.
Then for any element of $\Br(D)_3$ which is not in $\Br_0(D)$,
the evaluation map $D(\Q_3)\to\frac{1}{3}\Z/\Z$ is surjective.
\epr
{\em Proof}. By Lemma \ref{ex3} it is enough to assume that $D$ is one of
the surfaces $A$, $B$, $C$, and then assume that our element is 
$\sA\in \Br(D)_3$ constructed in Section \ref{calc}. 
Each of our surfaces has a $\Q$-point,
and at any such point the value of $\sA$ is zero by Lemma \ref{global}.
For each surface we shall exhibit a point 
$L=(x,y,z,w)$ in $D(\Q_5)$, where $x,y,z,w$ are integers without 
a common factor, such that neither $S_3(a,b:x,y)$ nor 
$S_3(c,d:z,w)$ holds, and use Corollary \ref{co2}.
The surjectivity of $\ev_{\sA,3}$ then follows from Corollary \ref{minus}.

Using Tate's algorithm one checks that $E^m$ has additive reduction
if and only if $4$ does not divide $\val_3(m)$. 

In the case of good reduction the group
of points on the reduced curve has no 3-torsion, hence
in this case we have $E^m(\Q_3)/3=E^m(\Q_3)_1/E^m(\Q_3)_2$.

\medskip

$A$. {\em In this case we prove more, namely, if $9$ does not divide $xyw$,
then neither $S_3(a,b:x,y)$ nor $S_3(c,d:z,w)$ holds.}
An example of such a point is $(\alpha,3,1,1)$ for an appropriate $\alpha\in\Z_3^*$.

Indeed, suppose that $\val_3(xyw)\leq 1$. Then an easy argument shows 
that $z$ is a unit and so is exactly one of $x$ and $y$.
Thus $\val_3(xy)=1$ and hence $w\in\Z_3^*$. Since $f\in\Z_3^*$
the curve $E'$ with equation $Y^2=X^3-4abf^2X$ has good reduction.
The point $P$ is in $E^m(\Q_3)_1$ but not in $E^m(\Q_3)_2$, so
the property $S_3(a,b:x,y)$ does not hold. The curve
$E''$ with equation $Y^2=X^3-4cdf^2X$ has additive reduction.
The point $Q+(0,0)$ with $X$-coordinate
$f^2/z^2w^2$ is in $E''(\Q_3)_0$ but not in $E''(\Q_3)_1$.
Because of (\ref{neron}) this point is not in $3E''(\Q_3)$.
Since $3(0,0)=(0,0)$, we conclude that $Q\notin 3E''(\Q_3)$.

\medskip

$B$. {\em In this case we prove more, namely, 
if $3$ does not divide $zw$,
then neither $S_3(a,b:x,y)$ nor $S_3(c,d:z,w)$ holds.}
An example of such a point is $(\alpha,1,1,1)$ for an appropriate $\alpha\in\Z_3^*$.

Indeed, suppose that $z,w\in\Z_3^*$.
Then $x,y\in\Z_3^*$ and $\val_3(f)=1$.
The curve $E'$ with equation $Y^2=X^3-4abf^2X$ has additive reduction.
The $X$-coordinate of $P$ is $-4abx^2y^2$, so $P$
is in $E'(\Q_3)_0$ but not in $E'(\Q_3)_1$, thus
$S_3(a,b:x,y)$ does not hold.

The curve $E''$ with equation $Y^2=X^3-4cdf^2X$ has additive reduction.
To get the minimal model we write $X=9X'$, $Y=27Y'$, then
the point $Q$ has $X'$-coordinate $\pm 12z^2w^2$. 
The point $Q+(0,0)$ has $X'$-coordinate $\pm f^2/9z^2w^2$,
and so $Q+(0,0)$ is in $E''(\Q_3)_0$ but not in $E''(\Q_3)_1$. 
Thus $Q+(0,0)$, and hence $Q$, is not in $3E''(\Q_3)$, so
that $S_3(c,d:z,w)$ does not hold. 

\medskip

$C$. Take $x=w=1$ and $y=2$.
Then there exists $z\in\Z_3^*$ such that $(x,y,z,w)\in D(\Q_3)$.
We have $f=17$, so the curve
$E'$ with equation $Y^2=X^3-(34)^2 X$ has good reduction.
The point $P$ has $X$-coordinate $-16$, so it reduces to the same point
as the 2-torsion point $M=(-34,0)$. Using {\tt sage} we check that
the valuation of the $X$-coordinate of $P+M$ is $-2$, hence
$P+M$ is in $E'(\Q_3)_1$ but not in $E'(\Q_3)_2=3E'(\Q_3)_1$.
The number of $\F_3$-points on the reduction is 4, thus 
$$E'(\Q_3)/3=E'(\Q_3)_1/3=E'(\Q_3)_1/E'(\Q_3)_2.$$ 
Therefore, $S_3(a,b:x,y)$ does not hold.

Next, the curve $E''$ with equation $Y^2=X^3-4cdf^2X$ has additive
reduction. The point $Q+(0,0)$ with $X$-coordinate $f^2/z^2w^2$
is in $E''(\Q_3)_0$ but not in $E''(\Q_3)_1$. Therefore, this
point, and hence also $Q$, is not in $3E''(\Q_3)$. Thus $S_3(c,d:z,w)$
does not hold.
$\Box$

\medskip

The proof of Theorem \ref{main2} is now complete.

\bigskip

\noindent Department of Mathematics and Statistics,
University of Cyprus, P.O. Box 20537, 1678, Nicosia, CYPRUS 
\medskip

\noindent ieronymou.evis@ucy.ac.cy

\bigskip

\noindent Department of Mathematics, South Kensington Campus,
Imperial College London, SW7 2BZ England, U.K.

\smallskip

\noindent Institute for the Information Transmission Problems,
Russian Academy of Sciences, 19 Bolshoi Karetnyi,
Moscow, 127994 Russia
\medskip

\noindent a.skorobogatov@imperial.ac.uk

\end{document}